\documentclass[10pt,reqno]{article}
\usepackage{amsmath,amssymb,amsthm}
\usepackage{esint}
\usepackage{bm}
\usepackage{url}
\usepackage{subfigure}
\usepackage{graphicx,color}
\usepackage{algorithm}
\usepackage{algpseudocode}
\usepackage{algorithmicx}
\usepackage{hyperref}
\usepackage{float}
\usepackage{booktabs,multirow}
\usepackage{hhline}
\usepackage{setspace}

\usepackage{caption}
\usepackage{hyphsubst}
\usepackage{caption}

\usepackage{float}
\usepackage{filecontents}
\usepackage{hyperref}

\hypersetup{colorlinks,citecolor=blue,linkcolor=blue, urlcolor=blue}
\usepackage{pdflscape}
\usepackage{breakurl}

\DeclareMathAlphabet{\itbf}{OML}{cmm}{b}{it}

\newtheorem{thm}{Theorem}[section]

\newtheorem{rem}[thm]{Remark}

\numberwithin{equation}{section}

\newcommand{\email}[1]{\protect\href{mailto:#1}{#1}}

\newcommand{\pathfigures}{Figures/}
\graphicspath{{\pathfigures}}

\begin{document}

\title{Comments on \textquotedblleft A new computing approach for power signal modeling using fractional adaptive algorithms\textquotedblright
}

\author{
Abdul Wahab\footnotemark[1]\, \footnotemark[2]
\and
Shujaat Khan\footnotemark[3]
\and 
Farrukh Zeeshan Khan\footnotemark[4]
}
\maketitle
\renewcommand{\thefootnote}{\fnsymbol{footnote}}
\footnotetext[1]{Corresponding Author. E-mail address:  \email{abdul.wahab@sns.nust.edu.pk}.}
\footnotetext[2]{Department of Mathematics, School of Natural Sciences, National University of Sciences and Technology (NUST), Sector H-12, 44000, Islamabad, Pakistan (\email{abdul.wahab@sns.nust.edu.pk)}.}
\footnotetext[3]{Bio-Imaging, Signal Processing, and Learning Lab., Department of Bio and Brain Engineering, Korea Advanced Institute of Science and Technology, 291 Daehak-ro, Yuseong-gu, 34141, Daejeon, South Korea (\email{shujaat@kaist.ac.kr}).}
\footnotetext[4]{Department of Computer Science, University of Engineering and Technology Taxila, 47080, Taxila, Pakistan (\email{farrukh.zeeshan@uettaxila.edu.pk)}.}
\renewcommand{\thefootnote}{\arabic{footnote}}

\begin{abstract}

In this paper, some points to the convergence analysis performed in the paper [A new computing approach for power signal modeling using fractional adaptive algorithms, ISA Transactions 68 (2017) 189-202]  are presented. It is highlighted that the way the authors prove convergence, suffers lack of correct and valid mathematical justifications. 
\end{abstract}


\noindent {\footnotesize {\bf Key words.} Power signal; Signal modeling; Parameter estimation; Fractional adaptive algorithms; Nonlinear adaptive strategies.}

\section{Introduction}\label{S:Intro}

In \cite{paper}, the fractional adaptive signal processing algorithms are utilized for identification of parameters in power signals by means of variants of fractional least means square (LMS) algorithms. The main contribution of \cite{paper} is the mathematical and computational performance analysis of fractional variants of LMS presented in \cite{ref48,ref37} when solving the problem discussed in \cite{ref8, ref7}. The most important part of the study is the mathematical convergence analysis which is presented in \cite[Sec. 3.5]{paper}. However, there are some mathematical errors which are detrimental to the correctness of the entire framework. The main objective of this note is to list those errors and substantiate that the entire convergence analysis is mathematically invalid. 
\begin{rem}
The symbols, notations and equation numbers used in this comment are consistent with \cite{paper}.
\end{rem}

\section{Mathematical Errors}\label{S:Errors}
Let us first give some remarks regarding the design of the fractional algorithms used in \cite{paper}.

\subsection{Issues in Algorithm Designs}
\begin{enumerate}
\item It is assumed (refer to Page 90, paragraph after \cite[Eq. 12]{paper}) that ``the fractional derivative of a constant is zero''. This is, in fact, not true in general; consult, for example, \cite{Kilbas}. Also, it can be easily verified from \cite[Eq. 14]{paper}, by substituting $p=0$ that the derivative of $1$ of order $0<fr<1$ is $t^{-fr}{\Gamma(1)}/{\Gamma(1-fr)}\neq 0$.  

\item \cite[Eq. 13]{paper} is derived by taking fractional derivative of order $fr$ of objective function \cite[Eq. 7]{paper}. In doing so, one requires fractional chain rule that involves several terms, evidently different from \cite[Eq. 13]{paper}. This has already been reported in \cite{NODY}, along with the correct form. Refer to \cite{Tarasov, Kilbas} for details of the fractional chain rule. Nevertheless, \cite[Eq. 13]{paper} can be seen as an approximation. So, we do not discuss this any further. 

\item \cite[Eq. 16]{paper} is a vector equation containing fractional power $\mathbf{w}^{1-fr}(t)$ of a vector $\mathbf{w}(t)$. The equation  is taken from \cite{FLMS} wherein the \textit{fractional power of a vector} is defined component-wise. Whenever, there is a negative element in the vector $\mathbf{w}$, $\mathbf{w}^{1-fr}(t)$ will render a complex output and impede the algorithm to converge to a real sought value. See \cite{Bershad, WS, NODY} for detailed discussions on this.  In the case when $\mathbf{w}^{1-fr}(t)$ is not defined componentwise, its sense needs to be specified since fractional powers of vectors are not defined in mathematics.

\item To avoid appearance of the complex outputs, a modulus is introduced in \cite[Eq. 17]{paper}, i.e.,  $\mathbf{w}^{1-fr}(t)$ is replaced by $|\mathbf{w}|^{1-fr}(t)$ (which makes it a scalar). However, a componentwise vector multiplication operator $\circ$ is introduced in the equation, that does not make sense because now there is only one vector in the product $\mathbf{u}(t)\circ |\mathbf{w}|^{1-fr}(t)$. Same remarks are relevant to \cite[Eqs. 21, 23, 24, 25]{paper}.
\end{enumerate}

It is worthwhile mentioning that all the aforementioned errors can be taken as approximations in algorithmic design, although, they can be nasty when we talk about convergence analysis. Therefore, we do not make any claim here and, instead, use this information to establish our claims regarding convergence analysis.

\subsection{Issues in Convergence Analysis}

\begin{enumerate}
\item In the convergence analysis, \cite[Eq. 32]{paper} is considered that involves $|\mathbf{\hat{\theta}}|^{1-fr}$, however, in the subsequent equation \cite[Eq. 34]{paper} the magnitude disappears. Note that
$|\mathbf{\hat{\theta}}|^{1-fr}$ is a scalar whereas $\mathbf{\hat{\theta}}^{1-fr}$ is a vector. Therefore, \cite[Eq. 32]{paper} does not provide \cite[Eq. 34]{paper}. Moreover, vector $\mathbf{\hat{\theta}}^{1-fr}$ is added to a scalar $1$ in \cite[Eq. 34]{paper} that does not make sense mathematically. Moreover, \cite[Eq. 34]{paper}  it also does not make sense because there is a product of  a vector with a vector (that is a dyad) and the resultant is added again to a vector.  

\item The binomial theorem \cite[Eq. 39]{paper} is inappropriate for vectors. In fact, it is not defined in Mathematics. Simple reason is that the multiplication of vectors is not commutative! Moreover, how would you interpret  $\theta^k_{opt}\Delta\theta(t)^{n-k}$? What is its direction? Can we add it to $\theta^{n-k}_{opt}\Delta\theta(t)^{k}$? What would be the meaning and direction of the resultant?  Not only this, \cite[Eq. 39]{paper} is used with a fractional exponent and that with negative values of the components of $\theta(t)^{k}$ will render complex outputs.  More interesting is the fact that if $\theta$ is a vector, mathematically, $\theta^2$ is a 2-tensor (matrix) or a dyad and $\theta^3$ is a 3-tensor and so on. So the summation is adding scalar, vectors and matrices and $k$ - tensors for $k\geq 3$! One can argue that these powers are taken componentwise, even then, the application of the binomial theorem is incorrect.

\item In view of the above, \cite[Eq. 40]{paper} is mathematically invalid. The rest of the convergence analysis is based on this equation. \cite[Eq. 41]{paper} has similar issues.

\item It is claimed that $p-R\theta_{opt}=0$ after \cite[Eq. 42]{paper}. However, there is no guarantee that the fractional optimal solution will converge to the Wiener solution. In fact, it is established that the optimal solution of the fractional LMS algorithm is not identical to the Wiener solution  (see \cite{Minima}). 

\item The notation $I(t+1)$ is undefined. Apparently, it is $E[\Delta\theta(t)]$. 

\item A function $F(\Delta\theta(t), fr)$ is defined in \cite[Eq. 44]{paper}. It is not clear whether it is a scalar field, vector field or a matrix valued function. Moreover, important is the fact that why should such a function exist at first place?

\item In \cite[Eq. 46]{paper} function $F$ is subtracted from the matrix $R$ and the resultant is added in a scalar $1$. Recall that we do not know whether $F$ is a scalar, vector, or a matrix. The same happens in \cite[Eq. 48]{paper} wherein an inequality is established containing scalar extremes sandwiching the difference of a scalar and difference of the matrix $R$ and function $F$. 

\item In \cite[Eq. 49, 50]{paper}, matrices are divided. The division of matrices is not defined in mathematics!

\end{enumerate}

Based on these observations, it is clear that the entire convergence analysis suffers lack of correct and valid mathematical justifications.  

\section{Conclusions} 
In this note, we have provided details of mathematical mistakes in the convergence analysis presented in \cite{paper}. Based on our observations, it is concluded that the convergence analysis in \cite{paper} is invalid and mathematically wrong.

\section*{Conflict of Interest} 
The authors declare that they have no conflict of interest.

\bibliographystyle{plain}

\end{document}